\documentclass[a4paper,12pt]{article}

\usepackage{listings}
\usepackage{longtable}

\usepackage{graphicx}

\usepackage[font=scriptsize]{caption}

\usepackage{float}

\usepackage{amsmath,amssymb,mathrsfs,eucal}

\usepackage[T1,T2A]{fontenc}

\font\cyrfont=wncyss10

\def\sza{\hbox{\cyrfont X}} 

\newtheorem{thm}{Theorem}
\newtheorem{cor}{Corollary}
\newtheorem{lem}{Lemma}
\newtheorem{prop}{Proposition}

\newtheorem{conj}[thm]{Conjecture}

\newtheorem{defin}[thm]{Definition}

\begin{document}

\title{Orders of Tate-Shafarevich groups for the 
Neumann-Setzer type elliptic curves}

\author{Andrzej D\k{a}browski and Lucjan Szymaszkiewicz} 

\date{}

\maketitle{}

{\it Abstract}. We present the results of our search for the 
orders of Tate-Shafarevich groups for the 
Neumann-Setzer type elliptic curves. 

\bigskip 
Key words: elliptic curves, Tate-Shafarevich group, 
Cohen-Lenstra heuristics, distribution 
of central $L$-values

\bigskip 
2010 Mathematics Subject Classification: 11G05, 11G40, 11Y50

\section{Introduction}

Let $E$ be an elliptic curve defined over $\Bbb Q$ of conductor $N_E$, 
and let $L(E,s)$ denote its $L$-series. 
Let $\sza(E)$ be the Tate-Shafarevich group of $E$, 
$E(\Bbb Q)$ the group of rational points, and $R(E)$ the regulator,  
with respect to the N\'eron-Tate height pairing. 
Finally, let $\Omega_E$ be the least positive real period of the N\'eron 
differential on $E$, and define $C_{\infty}(E)=\Omega_E$ or $2\Omega_E$ 
according as $E(\Bbb R)$ is connected or not, and let $C_{\text{fin}}(E)$ 
denote the product of the Tamagawa factors of $E$ at the bad primes. 
The Euler product defining $L(E,s)$ converges for $\text{Re}\,s>3/2$. 
The modularity conjecture, proven by Wiles-Taylor-Diamond-Breuil-Conrad,   
implies that $L(E,s)$ has an analytic continuation to an entire function. 
The Birch and Swinnerton-Dyer conjecture relates the arithmetic data 
of $E$ to the behaviour of $L(E,s)$ at $s=1$. 

Let $g_E$ be the rank of $E(\Bbb Q)$ and let $r_E$ denote the order of the zero 
of $L(E,s)$ at $s=1$. 

\begin{conj} (Birch and Swinnerton-Dyer) (i) We have $r_E=g_E$, 

(ii) the group $\sza(E)$ is finite, and 
$$ 
\lim_{s\to1}\frac{L(E,s)}{(s-1)^{r_E}} = 
\frac{C_{\infty}\, (E)C_{\text{fin}}(E)\, R(E)\, |\sza(E)|}{|E(\Bbb Q)_{\text{tors}}|^2}. 
$$
\end{conj} 
If $\sza(E)$ is finite, the work of Cassels and Tate shows that its order 
must be a square. 

The first general result in the direction of this conjecture was proven for 
elliptic curves $E$ with complex multiplication by Coates and Wiles in 1976 
\cite{CW}, who showed that if $L(E,1)\not =0$, then the group $E(\Bbb Q)$ is finite. 
Gross and Zagier \cite{GZ} showed that if $L(E,s)$ has a first-order zero at $s=1$, 
then $E$ has a rational point of infinite order. Rubin \cite{Rub} proves that 
if $E$ has complex multiplication and $L(E,1)\not =0$, then $\sza(E)$ is finite. 
Kolyvagin \cite{Kol} proved that, if $r_E\leq 1$, 
then $r_E=g_E$ and $\sza(E)$ is finite.  Very recently, Bhargava, Skinner 
and Zhang \cite{BhSkZ} proved that at least $66.48 \%$ of all elliptic curves 
over $\Bbb Q$, when ordered by height, satisfy the weak form of the Birch 
and Swinnerton-Dyer conjecture, and have finite Tate-Shafarevich group. 

Coates et al. \cite{cltz} \cite{coa}, and Gonzalez-Avil\'es 
\cite{G-A} showed that there is 
a large class of explicit quadratic twists of $X_0(49)$ whose 
complex $L$-series does not vanish at $s=1$, and for which the full 
Birch and Swinnerton-Dyer conjecture is valid. The deep results 
by Skinner-Urban \cite{SkUr} allow (in practice, see section 3 for  instance) 
to establish the full version of the Birch and Swinnerton-Dyer conjecture for 
a large class of elliptic curves without CM.

The numerical studies and conjectures by 
Conrey-Keating-Rubinstein-Snaith \cite{CKRS}, 
Delaunay \cite{Del0}\cite{Del}, Watkins \cite{Wat}, Radziwi\l\l-Soun\-dararajan \cite{RS}  
(see also the papers \cite{dw} \cite{djs} \cite{DSz} and references therein) 
substantially extend the systematic tables given by Cremona.

Given an integer $u\equiv 1 (\text{mod} \,4)$, such that 
$u^2+64$ is square-free, we define two  
families of elliptic curves of conductor 
$u^2+64$ (we call them  
the {\it Neumann-Setzer type elliptic curves}): 
$$
E_1(u): \quad y^2+xy=x^3+{1\over 4}(u-1)x^2-x, 
$$ 
and 
$$
E_2(u): \quad y^2+xy=x^3+{1\over 4}(u-1)x^2+4x+u. 
$$
In this paper we present the results of our search for the 
orders of Tate-Shafarevich groups for the 
Neumann-Setzer type elliptic curves. 
Our data contains values of $|\sza(E_i(u))|$ for $2056445$ values of 
$u\equiv 1 (\text{mod}\, 4)$, $|u|\leq 10^7$ such that $u^2+64$ 
is a product of odd number of different primes, 
and such that $L(E(u),1)\not=0$ 
($456702$ of these values satisfy the condition $u^2+64$ is a prime). 
Additionally, we have considered $10000$ values of 
$u\equiv 1 (\text{mod}\, 4)$, $|u|\geq 10^8$ such that $u^2+64$ 
is a product of odd number of different primes, 
and in cases $L(E(u),1)\not=0$ we computed the orders of 
$\sza(E_i(u))$. Our data extends the calculations given 
by Stein-Watkins \cite{SW} (resp. by Delaunay-Wuthrich \cite{DWut}), 
where the authors considered $|u|\leq \sqrt{2} \times 10^6$ 
(resp. $|u|\leq 10^6$) such that $u^2+64$ is a prime. 

Our main observations concern the asymptotic 
formulae in sections 4 (frequency of orders of $\sza$) and 6 
(asymptotics for the sums $\sum |\sza(E_i(u))| R(E_i(u)$ in the 
rank zero and one cases), and the distributions of 
$\log L(E_i(u),1)$ and $\log (|\sza(E_i(u))|/\sqrt{|u|})$ in 
section 7. 
 
\bigskip 

We thank Bjorn Poonen and Christophe Delaunay for their remarks 
and questions. We thank the anonymous referee for his/her remarks 
and comments which improved the final version of this paper.

\bigskip 
Our experimental data were obtained using the the PARI/GP software 
\cite{PARI}. 
The computations were carried out in 2015 and 2016 
on the HPC cluster HAL9000 and 
desktop computers Core(TM) 2 Quad Q8300 4GB/8GB. 
All machines are located at the Department of Mathematics and 
Physics of Szczecin University.

\section{Preliminaries}

We have $\Delta_{E_1(u)}=u^2+64$, and  $\Delta_{E_2(u)}=-(u^2+64)^2$. 
The curves $E_1(u)$ and $E_2(u)$ are $2$-isogenous: write $E_1(u)$ and 
$E_2(u)$ in short Weierstrass forms ($y^2=x^3+ux^2-16x$ and 
$y^2=x^3-2ux^2+(u^2+64)x$, respectively), and use (\cite{Sil}, 
Example 4.5 on p. 70).  
It is known, due to Neumann and Setzer (\cite{Neu}, \cite{Set}), 
that in the case $u^2+64$ is a prime, the curves 
$E_1(u)$ and $E_2(u)$ are the only (up to isomorphism) 
elliptic curves 
with a rational $2$-division point and conductor $u^2+64$. 
In general there are more than two, up to isomorphism, 
elliptic curves with a rational $2$-division point and 
conductor $u^2+64$. Take, for instance, $u=-51$, then 
the curves $E_1(u)$ and $E_2(u)$ have conductor 
$2665=5\cdot 13\cdot 41$. In Cremona's online tables 
we find $8$ elliptic curves of conductor 
$2665$ with a rational $2$-division point.

\begin{lem} We have 
(i) $E_1(u)(\Bbb Q)_{tors}\simeq E_2(u)(\Bbb Q)_{tors} 
\simeq \Bbb Z/2\Bbb Z$; 
(ii) $\Omega_{E_1(u)}=\Omega_{E_2(u)}$,  
$C_{\infty}(E_1(u))=2\Omega_{E_1(u)}$, $C_{\infty}(E_2(u))=\Omega_{E_2(u)}$; 
(iii) $C_{fin}(E_1(u))=1$, and $C_{fin}(E_2(u))=2^k$, 
where $u^2+64=p_1\cdots p_k$. 
\end{lem} 

\noindent
{\it Proof}. (i) Let $E(u)=E_1(u)$ or $E_2(u)$. Then $E(u)$ 
has good reduction at $2$. 
Using the reduction map modulo $2$, we obtain that 
$|E_i(u)(\Bbb Q)_{tors}|$ divides $4$. Now, one checks that 
$E_i(u)(\Bbb Q)$ 
have only one point of order two, and  no points of order four. 
(ii) To check that $\Omega_{E_1(u)}=\Omega_{E_2(u)}$, 
one uses the explicit forms of Weierstrass equations. 
Now the sign of the discriminant of $E_1(u)$ (resp. of $E_2(u)$) 
is positive (resp. negative), hence the remaining assertions follow. 
(iii) We have 
$C_{fin}(E_1(u))=\prod_{p|\Delta_{E(u)}}C_p(E(u))$, 
where $C_p(E(u))=[E(u)(\Bbb Q_p):E_0(u)(\Bbb Q_p)]$, and  
$E_0(u)(\Bbb Q_p)$ denotes the subgroup of points of 
$E(u)(\Bbb Q_p)$ with non-singular reduction modulo $p$. 
Both $E_1(u)$ and $E_2(u)$ have split 
multiplicative reductions at all primes $p$ dividing 
$u^2+64$. Hence, in this case, 
$C_p(E(u))=\text{ord}_p(\Delta_{E(u)})$ (see, for instance, 
\cite{coa}, Lemma 2.9), and the assertion follows. 

\bigskip 
Note that $L(E_1(u),s) =  L(E_2(u),s) =  \sum_{n=1}^{\infty}a_nn^{-s}$, 
$\text{Re}(s) > 3/2$. 
Assuming the truth of the Birch and Swinnerton-Dyer 
conjecture for $E(u)$ in the rank zero case, we can 
calculate the order of $\sza (E(u))$ by evaluating 
(an analytic continuation of) $L(E(u),s)$ at $s=1$: 
$$
|\sza(E_1(u))| = {2 L(E_1(u),1) \over\Omega_{E_1(u)}},  
$$
$$
|\sza(E_2(u))| = {L(E_2(u),1) \over 2^{k-2}\Omega_{E_2(u)}},  
$$
where as above, $u^2+64=p_1\cdots p_k$ is a product of 
different primes. 

More precisely, we have to calculate the value 
$$
L(E(u),1) = 2\sum_{n=1}^{\infty}{a_n\over n}\, 
e^{-{2\pi n\over \sqrt{u^2+64}}} 
$$
with sufficiently accuracy. 

\begin{lem} In order to determine the order of $\sza(E_1(u))$ 
and $\sza(E_2(u))$, 
it is enough to take ${1\over 8}\sqrt{u^2+64}\log(u^2+64)$ 
terms of the above series. 
\end{lem} 

\noindent
{\it Proof}. Repeat the proof of Theorem 16 in \cite{DWut}.

\bigskip 
Let $\epsilon(E(u))$ denote the root number of $E(u)$.  

\begin{lem} Let $u^2+64=p_1\cdots p_k$ be a product of 
different primes. Then $\epsilon(E(u))=(-1)^{k+1}$. 
\end{lem} 

\noindent
{\it Proof}. $\epsilon(E(u))=
-\prod_{i=1}^k\epsilon_{p_i}(E(u))$, a product of local 
root numbers. Now, $E(u)$ has split multiplicative reduction 
at all $p_i$ dividing $u^2+64$. Hence, $\epsilon_{p_i}(E(u))=-1$, 
and the assertion follows.

\begin{cor} Assume the parity conjecture holds for the 
curves $E(u)$. Then $E(u)(\Bbb Q)$ has even rank if and only 
if  $u^2+64=p_1\cdots p_k$ is a product of odd number of 
different primes. 
\end{cor}

We can use a classical 2-descent method (\cite{Sil}, Chapter X) 
to obtain a bound on the rank of $E_i(u)$ depending on $k$. Let 
$\phi: E_1(u) \to E_2(u)$ be the $2$-isogeny, and write $\hat\phi$ 
for its dual. Let $S^{(\phi)}$ and $S^{(\hat\phi)}$ denote the 
corresponding Selmer groups. One checks that 
$S^{(\phi)}\subset <p_1,...,p_k>$ and $S^{(\hat\phi)}=<-1>$. 
As a consequence, we obtain 
$\text{rank}(E_i(u)) \leq \dim_{\Bbb F_2}S^{(\phi)} + 
\dim_{\Bbb F_2}S^{(\hat\phi)} -2 \leq k + 1 - 2 = k-1$. 
In particular, if $u^2+64$ is a prime, then $E_i(u)$ have rank zero, 
and if $k=2$, then  $\text{rank}(E_i(u)) \leq 1$ ($=1$ if we assume 
the parity conjecture).

\begin{defin} We say that an integer 
$u\equiv 1 (\text{mod}\, 4)$ 
satisfies condition (*), if $u^2+64$ is a prime; we say that 
an integer $u\equiv 1 (\text{mod}\, 4)$ 
satisfies condition (**), if $u^2+64$ is a product of 
odd number of different primes. 
\end{defin}

\section{Birch and Swinnerton-Dyer conjectur for Ne\-umann-Setzer 
type elliptic curves}

In this section, we will use the deep results by Skinner-Urban 
\cite{SkUr} (and other available techniques), to prove the full 
version of the Birch-Swinnerton-Dyer conjecture for a large 
class of Neumann-Setzer type elliptic curves.  

Let $\overline{\rho}_{E,p}: 
\text{Gal}(\overline{\Bbb Q}/\Bbb Q) \to 
\text{GL}_2(\Bbb F_p)$ denote the Galois 
representation on the $p$-torsion of $E$. 
Assume $p\geq 3$. 

\begin{thm} (\cite{SkUr}, Theorem 2) Let $E$ be an 
elliptic curve over $\Bbb Q$ with conductor $N_E$. 
Suppose: (i) $E$ has good ordinary reduction at $p$; 
(ii) $\overline{\rho}_{E,p}$ is irreducible; (iii) 
there exists a prime $q\not=p$ such that $q\mid\mid N_E$ 
and $\overline{\rho}_{E,p}$ is ramified at $q$; 
(iv) $\overline{\rho}_{E,p}$ is surjective.  
If moreover $L(E,1)\not=0$, then the $p$-part of the Birch 
and Swinnerton-Dyer conjecture holds true, and  
we have 
$$
\text{ord}_p(|\sza(E)|)=\text{ord}_p\left(
{|E(\Bbb Q)_{\text{tors}}|^2 L(E,1)\over C_{\infty}(E)C_{fin}(E)}\right). 
$$
\end{thm}

Take $E(u)=E_1(u)$ or $E_2(u)$. Then: 

a) $E(u)$ is semistable and has a rational $2$-division 
point, hence $\overline{\rho}_{E(u),p}$ is irreducible for 
$p\geq 7$ by (\cite{DM}, Theorem 7). Note moreover 
(by Wiles \cite{Wi}) that at least one of 
$\overline{\rho}_{E(u),3}$ or $\overline{\rho}_{E(u),5}$ is 
irreducible. 

b) If $E$ is any semistable elliptic curve and $q\not=p$, 
then $\overline{\rho}_{E,p}$ is unramified at $q$ if 
and only if  
$p|\text{ord}_q(\Delta_E)$. In our case, 
$\text{ord}_q(\Delta_E(u))$ equals $1$ or $2$, hence 
$\overline{\rho}_{E(u),p}$ is ramified at any $q\geq 3$.  

c) If $E$ is any semistable elliptic curve, then 
$\overline{\rho}_{E,p}$ is surjective for $p\geq 11$ 
by \cite{Se}. More precisely, Serre (\cite{Se}, Prop. 1) 
shows that in this case $\overline{\rho}_{E,p}$ is surjective 
for all primes $p$ unless $E$ admits an isogeny of degree 
$p$ defined over $\Bbb Q$. In particular, if such $E$ 
additionally has a rational $2$-division point, then 
$\overline{\rho}_{E,p}$ is surjective for $p\geq 7$. 
Note (by \cite{Se0}, Prop. 21, and \cite{Se}, Prop. 1), 
that in the case of semistable elliptic curve $E$,  
the representation $\overline{\rho}_{E,p}$ is surjective 
if and only if it is irreducible. Now, Zywina (\cite{Zyw}, Prop. 6.1) 
gives a criterion to determine whether $\overline{\rho}_{E,p}$ 
is surjective or not for any non-CM elliptic curve and any prime 
$p\leq 11$. Using such a criterion, one immediately checks 
surjectivity of $\overline{\rho}_{E_i(u),p}$ for $p=2$, $3$, and $5$. 
As a consequence, we obtain the following general result.

\begin{prop}
The representations $\overline{\rho}_{E(u),p}$ 
are surjective for all primes $p$.   
\end{prop} 

Summing up all the above information, we obtain the following 
nice result. 

\begin{cor} 
Let $E=E_1(u)$ or $E_2(u)$, 
with $u\equiv 1 (\text{mod}\, 4)$ satisfying (**) 
and such that $L(E,1)\not=0$. If $E$ has good ordinary reduction 
at $p\geq 3$, then the $p$-part of the Birch and Swinnerton-Dyer  
conjecture holds for $E$. 
\end{cor} 

\noindent 
{\bf Remark.} Let us recall that a prime $p$ is {\it good} for an elliptic 
curve $E$ over $\Bbb Q$, if $p$ does not divide $N_E$; $p$ is {\it good 
ordinary} for $E$, if is good and $a_p=p+1-N_p(E)$ is not divisible 
by $p$ (here $N_p(E)$ denotes the number of $\Bbb F_p$-points of the 
reduction $E_p$). Here are explicit conditions for small primes $p$ 
to satisfy the good ordinary condition in case $E=E_i(u)$ 
(we assume $u\equiv 1 (\text{mod} \,4)$): 
(i) $p=3$, additional condition 
$u\not\equiv 0 (\text{mod} \,3)$; (ii) $p=5$, no additional condition on $u$; 
(iii) $p=7$, aditional condition 
$u\not\equiv 0 (\text{mod} \,7)$; 
(iv) $p=11$, additional condition $u\not\equiv 0, 4, 7 (\text{mod} \,11)$.

\bigskip 
\noindent 
{\bf Remark.} One can use explicit descent algorithms to compute 
$\sza(E_i(u))[m]$ for $m=2$, $4$ or $8$. If $\sza(E_i(u))[2]$ is 
trivial, then $\sza(E_i(u))$ has odd order. If 
$\sza(E_i(u))[2] = \sza(E_i(u))[4]$, say, then 
$\text{ord}_2|\sza(E_i(u))| = \text{ord}_2|\sza(E_i(u))[2]|$. 
Similarly, one can use explicit descent algorithms to compute 
$\sza(E_i(u))[m]$ for $m=3$ or $9$. Again, if $\sza(E_i(u))[3]$ 
is trivial, then $\sza(E_i(u))$ has order not divisible by $3$ 
(here we not require that $3$ is good ordinary). 
If $\sza(E_i(u))[3] = \sza(E_i(u))[9]$, then 
$\text{ord}_3|\sza(E_i(u))| = \text{ord}_3|\sza(E_i(u))[3]|$. 

The theses \cite{Mil} \cite{Soh} explore both theoretical and computational methods 
to compute the orders of Tate-Shafarevich groups.

\bigskip 
\noindent 
{\bf Remark.} (i) Among $456702$ values of 
$u\equiv 1 (\text{mod}\, 4)$, 
$|u|\leq 10^7$ satisfying (*), there are $379898$ values 
of $|u|$ such that $E(u)$ has good ordinary reduction at 
any prime dividing the analytic order $|\sza(E(u))|$. 
The groups $\sza(E_i(u))[2]$ are both trivial (by $2$-descent), 
hence by Corollary 2 the values $|\sza(E(u))|$ are the algebraic orders 
of $\sza$.  (ii) Among  $2056445$ values of $u\equiv 1 (\text{mod}\, 4)$, 
$|u|\leq 10^7$ satisfying (**) and such that $L(E(u),1)\not=0$, 
there are $1148683$ values 
of $|u|$ such that $|\sza(E_2(u))|$ is odd and $E(u)$ has good 
ordinary reduction at any prime dividing the analytic order 
$|\sza(E_2(u))|$. Again, by Corollary 2 all these values are the 
algebraic orders of $\sza$.

\bigskip 

The numerical data are done under the Birch and Swinnerton-Dyer conjecture. 
In particular, the experimental study in sections 4, 5, 6, and 7 
concern the analytic orders of the Tate-Shafarevich groups.

\section{Frequency of orders of $\sza$}

Our calculations strongly suggest that 
for any positive integer $k$ there are infinitely many 
integers $u\equiv 1 (\text{mod}\, 4)$ satisfying condition 
(**), such that $E(u)$ has rank zero and $|\sza(E(u))| = k^2$. 
Below (end of this section) we will state a more precise 
conjecture.

Let $f(i,X)$ denote the number of integers 
$u\equiv 1 (\text{mod}\, 4)$, $|u|\leq X$, satisfying (**) and 
such that $L(E(u),1)\not=0$, $|\sza(E_i(u))| = 1$. Let 
$g(X)$ denote the number of integers $u\equiv 1 (\text{mod}\, 4)$, $|u|\leq X$, 
satisfying (**) and such that $L(E(u),1) = 0$. 
We obtain the following graphs (compare \cite{djs} \cite{DSz}, 
where similar observations are made for the families of 
quadratic twists of several elliptic curves).

\begin{figure}[H]  
\centering
\includegraphics[trim = 0mm 20mm 0mm 15mm, clip, scale=0.4]{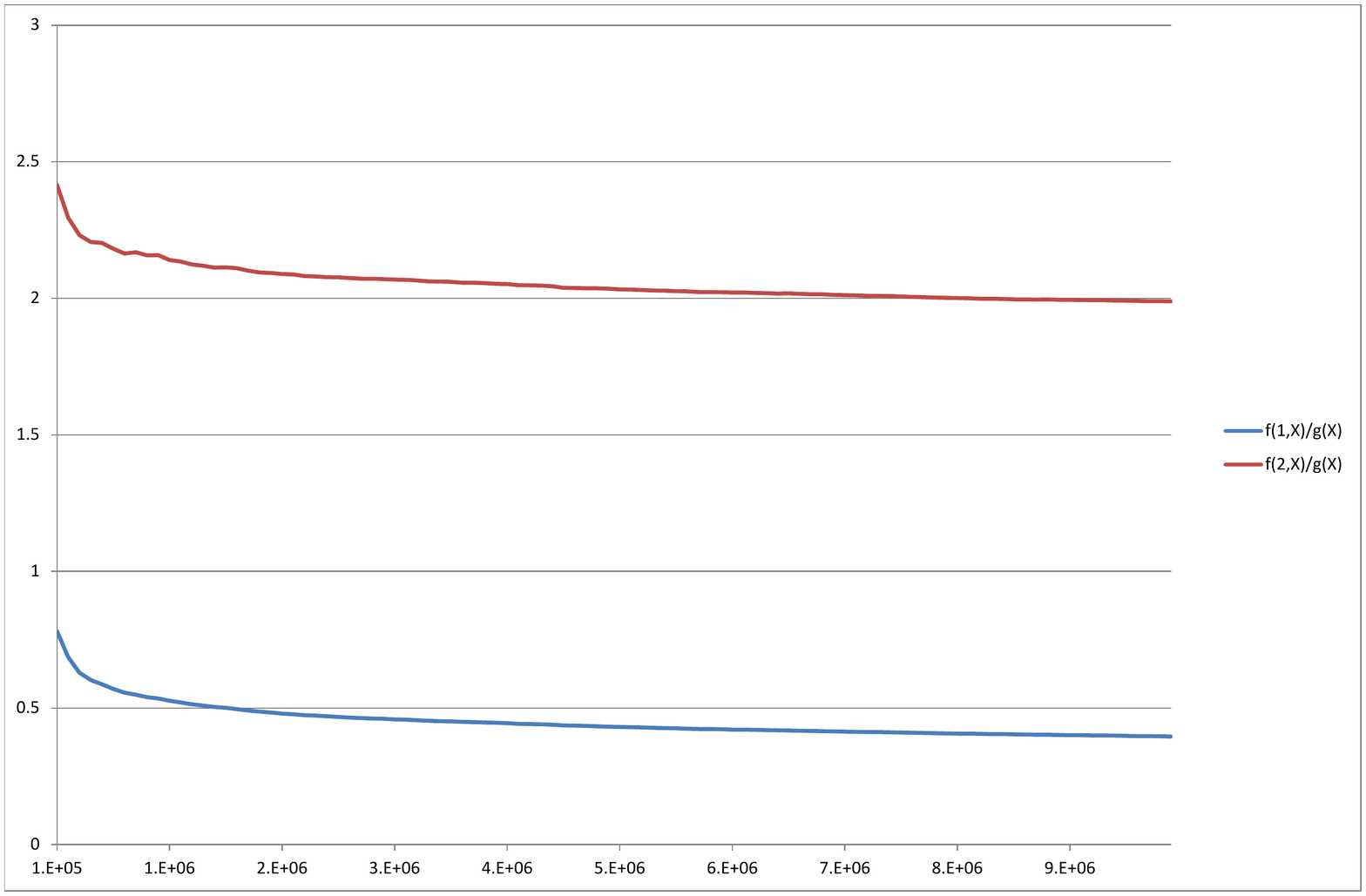}
\caption{Graphs of the functions $f(i,X)/g(X)$, $i=1,2$.} 
\end{figure}

Consider the set consisting of $10000$ values of integers 
$u\equiv 1 (\text{mod}\, 4)$, $|u|\geq 10^8$, satisfying (**). 
Let $f(i)$ denote the number of such $u$'s satisfying 
$L(E_i(u),1)\not=0$ and $|\sza(E_i(u))|=1$, and let $g$ 
denote the number of such $u$'s satisfying $L(E_i(u),1)=0$. 
Then $f(1)=118$, $f(2)=845$, $g=482$, hence $f(1)/g\approx 0,2448$, 
and $f(2)/g\approx 1,7531$.

Delaunay and Watkins expect \cite{DW}, Heuristics 1.1):  
$$
\sharp\{d\leq X: \epsilon(E_d)=1, \,\text{rank}(E_d)\geq 2\} 
\sim c_E X^{3/4}(\log X)^{b_E+{3\over 8}}, \quad \text{as} 
\quad X\to\infty, 
$$
where $c_E>0$, and there are four different possibilities for $b_E$, 
largery dependent on the rational $2$-torsion structure of $E$. 
Watkins \cite{Wat}, and Park-Poonen-Voight-Wood \cite{PPVW}  
have conjectured that 
$$
\sharp\{E: \text{ht}(E)\leq X, \, \epsilon(E)=1, \,\text{rank}(E)\geq 2\} 
\sim c X^{19/24}(\log X)^{3/8}, 
$$
where $E$ runs over all elliptic curves defined over the rationals, 
and $\text{ht}(E)$ denotes the height of $E$. 

We expect a similar asymptotic formula for the family $E(u)$. 
Let $H(X):={X^{19/24}(\log X)^{3/8}\over g(X)}$, and 
$G_i(X):={X^{3/4}(\log X)^{i}\over g(X)}$, $i=0$, $1/2$ or $1$.  
We obtain the following graphs, (partially) confirming our expectation.

\begin{figure}[H]  
\centering
\includegraphics[trim = 0mm 20mm 0mm 15mm, clip, scale=0.4]{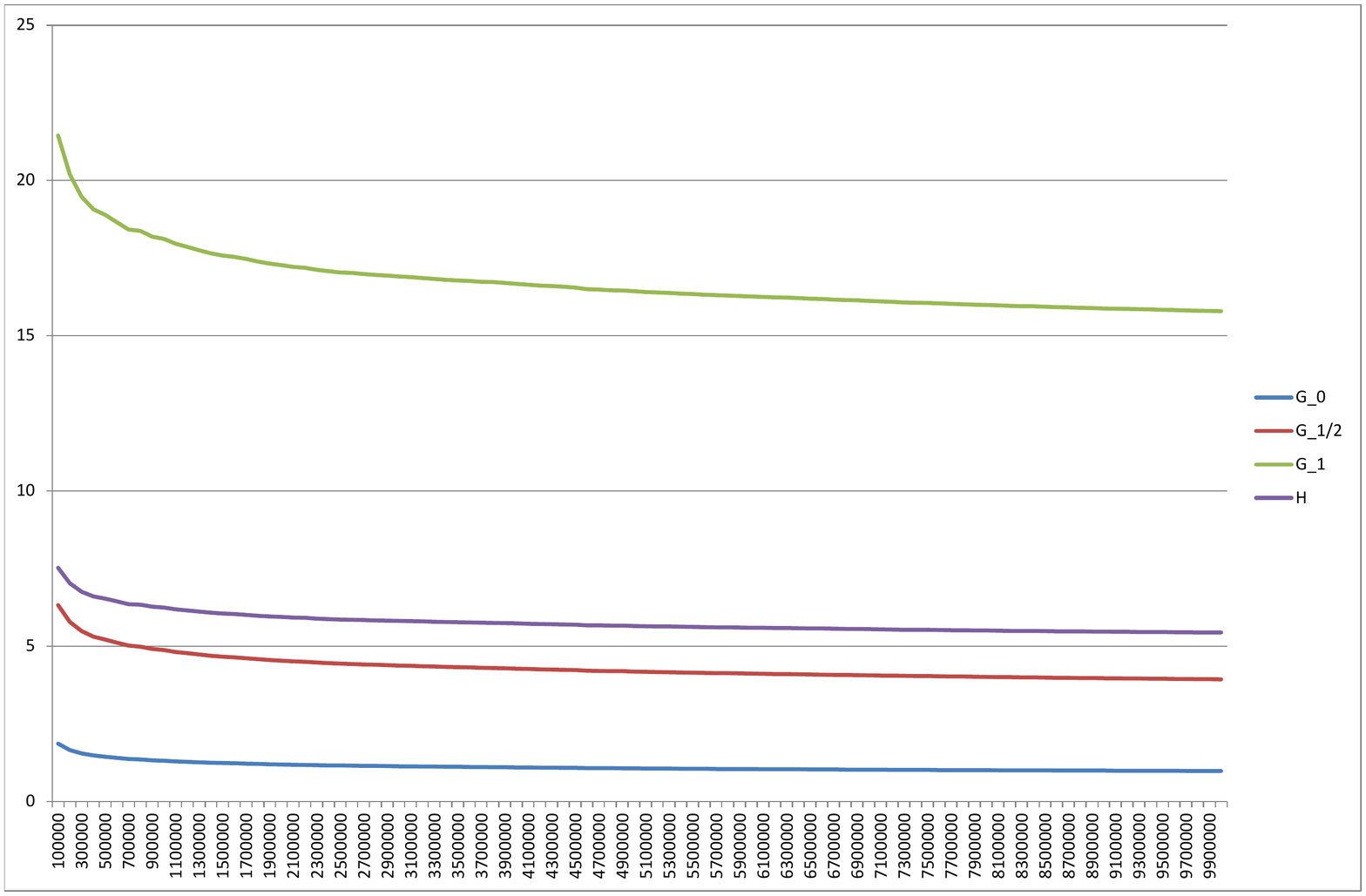}
\caption{Graph of the function $H(X)$.} 
\end{figure}

Now let $f_k(i,X)$ denote the number of integers 
$u\equiv 1 (\text{mod}\, 4)$, $|u|\leq X$, satisfying (**) and 
such that $L(E(u),1)\not=0$, $|\sza(E_i(u))| = k^2$. 
Let $F_k(i,X):={f(i,X)\over f_k(i,X)}$. We obtain the following graphs 
of the functions $F_k(i,X)$ for $i=1,2$ and $k=2,3,4,5,6,7$.

\begin{figure}[H]  
\centering
\includegraphics[trim = 0mm 20mm 0mm 15mm, clip, scale=0.4]{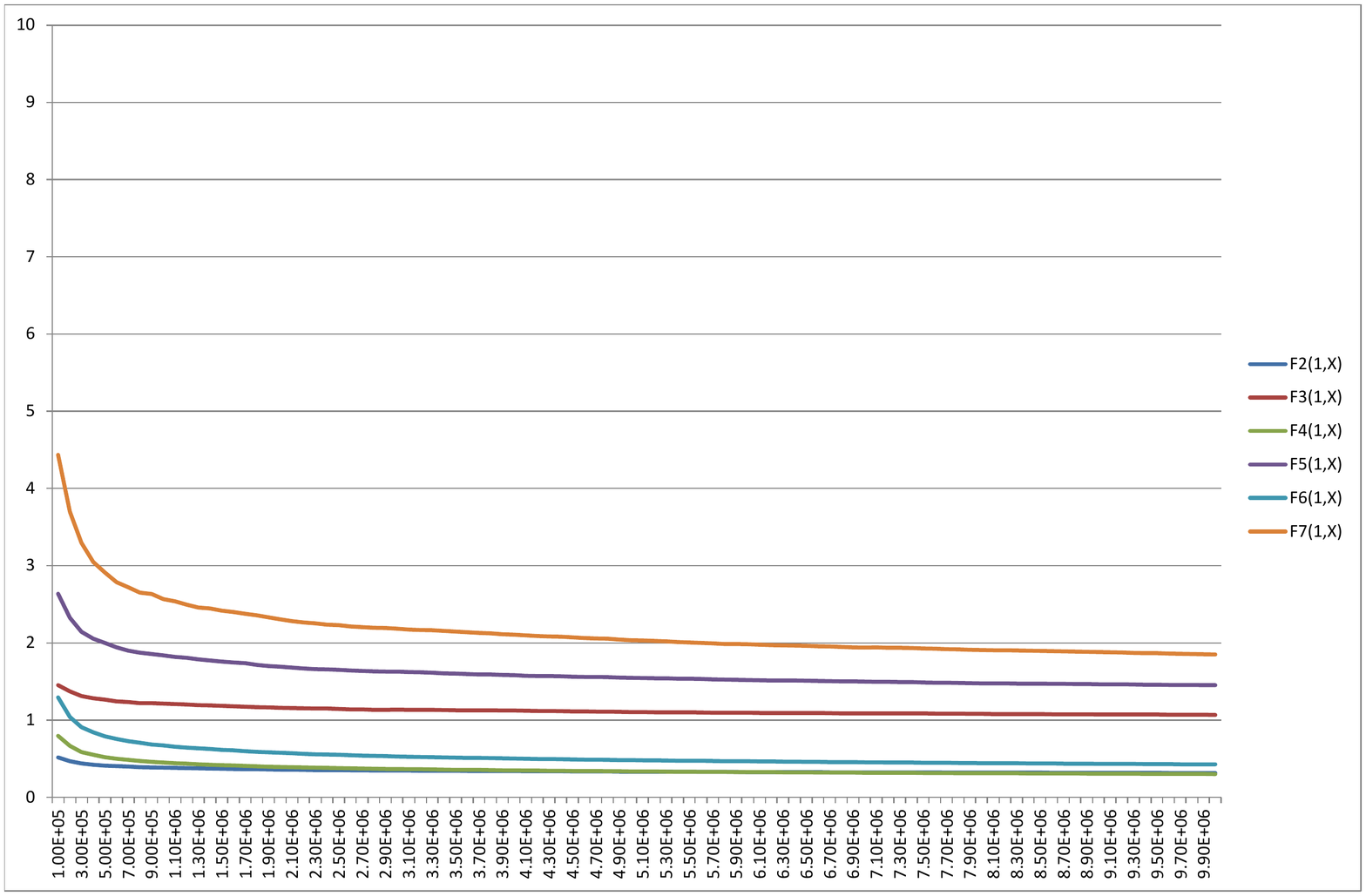}
\caption{Graphs of the functions $F_k(1,X)$, $k=2,...7$.} 
\end{figure} 

\begin{figure}[H]  
\centering
\includegraphics[trim = 0mm 20mm 0mm 15mm, clip, scale=0.4]{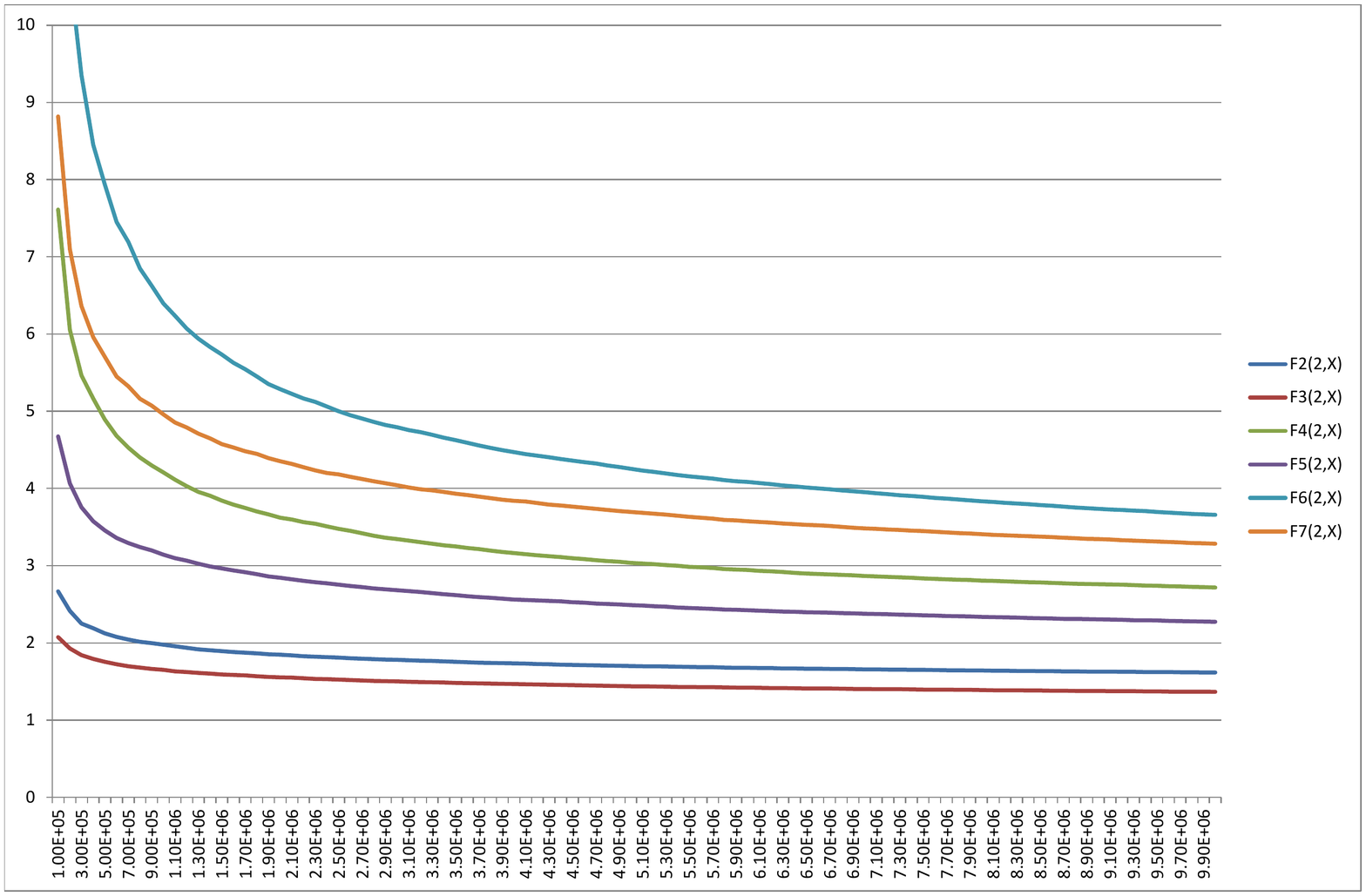}
\caption{Graphs of the functions $F_k(2,X)$, $k=2,...7$.} 
\end{figure}  

The above calculations suggest the following 

\begin{conj}
For any positive integer $k$ there are constants $c_{k,i}>0$, 
$\alpha_{k,i}$, and $\beta_{k,i}$ such that 
$$
f_k(i,X) \sim c_{k,i}X^{\alpha_{k,i}}(\log X)^{\beta_{k,i}}, \quad \text{as} 
\quad X\to\infty. 
$$
\end{conj}

Conjectures 8 in \cite{djs} and 2 in \cite{DSz} suggest similar asymptotics 
for the family of quadratic twists of any elliptic curve defined over $\Bbb Q$.

Consider the set consisting of $10000$ values of integers 
$u\equiv 1 (\text{mod}\, 4)$, $|u|\geq 10^8$, satisfying (**). 
Let $f_k(i)$ denote the number of such $u$'s satisfying 
$L(E_i(u),1)\not=0$ and $|\sza(E_i(u))|=k^2$.  
Let $F_k(i):={f_1(i)\over f_k(i)}$. We obtain 

\bigskip 

$F_2(1)\approx 0.2256$, \quad  $F_3(1)\approx 0.8251$, 
\quad $F_4(1)\approx 0.1779$  

$F_5(1)\approx 1.0825$, \quad $F_6(1)\approx 0.2494$, 
\quad  $F_7(1)\approx 1.1919$ 

$F_2(2)\approx 1.1901$, \quad  $F_3(2)\approx 1.0682$, 
\quad $F_4(2)\approx 1.5590$ 

$F_5(2)\approx 1.4955$, \quad $F_6(2)\approx 1.9031$, 
\quad  $F_7(2)\approx 1.8449$

\section{Cohen-Lenstra heuristics for the order of $\sza$} 

Delaunay \cite{Del} has considered Cohen-Lenstra heuristics for the order of 
Tate-Shafarevich group. He predicts, among others, that in the rank zero case, 
the probability that $|\sza(E)|$ of a given elliptic curve $E$ over $\Bbb Q$ is 
divisible by a prime $p$ should be 
$
f_0(p):=1-\prod_{j=1}^{\infty}(1-p^{1-2j})={1\over p}+{1\over p^3}+... \,.    
$
Hence, $f_0(2)\approx 0.580577$, $f_0(3)\approx 0.360995$, $f_0(5)\approx 0.206660$, 
$f_0(7)\approx 0.145408$,  $f_0(11)\approx 0.092$, 
and so on.

Let $F(X)$ (resp. $G(X)$) denote the number of integers 
$u\equiv 1 (\text{mod}\, 4)$, $|u|\leq X$,  
satisfying (*) (resp. (**)) and such that $L(E(u),1)\not=0$. 
Let $F_p(X)$ (resp. $G_p(X)$ if $p\geq 3$) denote 
the number of integers 
$u\equiv  1 (\text{mod}\, 4)$, $|u|\leq X$, satisfying (*) 
(resp. (**)), such that $L(E(u),1)\not=0$ and  
$|\sza(E(u))|$ is divisible by $p$. 
Let $G_2(i,X)$ denote the number of integers 
$u\equiv  1 (\text{mod}\, 4)$, $|u|\leq X$, satisfying (**), 
such that $L(E(u),1)\not=0$ and  
$|\sza(E_i(u))|$ is divisible by $2$. 
Let $f_p(X):={F_p(X)\over F(X)}$,  
$g_p(X):={G_p(X)\over G(X)}$, and $g_2(i,X):={G_2(i,X)\over G(X)}$. 
We obtain the following tables, extending the calculations 
given by Stein-Watkins \cite{SW}  and Delaunay-Wuthrich \cite{DWut}.

\begin{center}
\footnotesize
\begin{longtable}{|r|r|r|r|r|r|r|} 
\hline 
\multicolumn{1}{|c|}{$X$} 
& 
\multicolumn{1}{|c|}{$f_3(X)$} 
& 
\multicolumn{1}{|c|}{$f_5(X)$} 
& 
\multicolumn{1}{|c|}{$f_7(X)$} 
& 
\multicolumn{1}{|c|}{$f_{11}(X)$} 
\\ 
\hline 
\endhead 
\hline 
\multicolumn{5}{|r|}{{Continued on next page}} \\
\hline
\endfoot

\hline 
\hline
\endlastfoot 

$2 \cdot 10^6$ & 0.358355 &  0.189909 &  0.123182 &  0.061527 \\
$4 \cdot 10^6$ & 0.362001 &  0.192343 &  0.126864 &  0.066945 \\
$6 \cdot 10^6$ & 0.363294 &  0.194413 &  0.129213 &  0.069780 \\
$8 \cdot 10^6$ & 0.364051 &  0.196239 &  0.130556 &  0.071144 \\
$           10^7$ & 0.365067 &  0.197048 &  0.131812 &  0.072358 \\

\end{longtable}
\end{center}

The numerical values of $f_3(X)$ exceed the expected value $f_0(3)$. 
In general, the values $f_k(X)$ may tend to some constants depending 
on the various congruential values of $u$ (compare \cite{SW}). 

It seems that it would be better to consider $u$'s satisfying (**), 
but here the  convergence is very slow. Here are the results.

\begin{center}
\footnotesize
\begin{longtable}{|r|r|r|r|r|r|r|} 
\hline 
\multicolumn{1}{|c|}{$X$} 
& 
\multicolumn{1}{|c|}{$g_2(1,X)$} 
& 
\multicolumn{1}{|c|}{$g_2(2,X)$} 
& 
\multicolumn{1}{|c|}{$g_3(X)$} 
& 
\multicolumn{1}{|c|}{$g_5(X)$} 
& 
\multicolumn{1}{|c|}{$g_7(X)$} 
& 
\multicolumn{1}{|c|}{$g_{11}(X)$} 
\\ 
\hline 
\endhead 
\hline 
\multicolumn{7}{|r|}{{Continued on next page}} \\
\hline
\endfoot

\hline 
\hline
\endlastfoot 

$2 \cdot 10^6$ &  0.746231 & 0.313111 & 0.295592 &  0.127626 &  0.072959 &  0.030979 \\
$4 \cdot 10^6$ &  0.761104 & 0.326554 & 0.303529 &  0.134259 &  0.078513 &  0.034796 \\
$6 \cdot 10^6$ &  0.768805 & 0.333854 & 0.307670 &  0.138168 &  0.081543 &  0.036884 \\
$8 \cdot 10^6$ &  0.774040 & 0.338854 & 0.310603 &  0.140959 &  0.083638 &  0.038350 \\
$           10^7$ &  0.777917 & 0.342322 & 0.312758 &  0.143060 &  0.085332 &  0.039481 \\

\end{longtable}
\end{center} 

Note that the value $(g_2(1,10^7)+g_2(2,10^7))/2 \approx 0.56012$ is not so 
far from the expected one.  

\bigskip 

We have computed the orders of $9518$ pairs of Tate-Shafarevich 
groups $(\sza(E_1(u)), \sza(E_1(u)))$ for $|u|\geq 10^8$, 
$u\equiv 1 (\text{mod}\, 4)$, satisfying (**), 
and such that $L(E(u),1)\not=0$. We obtained the following table.

\begin{center}

\footnotesize

\begin{longtable}{|r|r|r|r|r|r|}

\hline

\multicolumn{1}{|c|}{$p$}

&

\multicolumn{1}{|c|}{$2$}

&

\multicolumn{1}{|c|}{$3$}

&

\multicolumn{1}{|c|}{$5$}

&

\multicolumn{1}{|c|}{$7$}

&

\multicolumn{1}{|c|}{$11$}

\\

\hline

\endhead

\hline

\multicolumn{6}{|r|}{{Continued on next page}} \\

\hline

\endfoot

\hline

\hline

\endlastfoot

Frequency of $p| \, |\sza(E_1(u))|$ & 0.826329 & 0.332213 & 0.167262 &
0.111053 & 0.058100 \\

\hline

Frequency of $p| \, |\sza(E_2(u))|$ & 0.393045 & 0.332213 & 0.167262 &
0.111053 & 0.058100 \\

\end{longtable}

\end{center}

\section{Asymptotic formulae}

\subsection{The rank zero case}

Let 
$M^*(T):={1\over T^*}\sum |\sza(E(u))|$, 
where the sum is over integers  $u\equiv 1 (\text{mod}\, 4)$, 
$|u|\leq T$, satisfying (*) and $L(E(u),1)\not=0$, 
and $T^*$ denotes the number of 
terms in the sum. Similarly, let 
$N_i^{**}(T):={1\over T_i^{**}}\sum |\sza(E_i(u))|$, 
where $i=1,2$, and the sum is over integers  
$u\equiv 1 (\text{mod}\, 4)$, $|u|\leq T$, satisfying (**) 
and $L(E(u),1)\not=0$, 
and $T_i^{**}$ denotes the number of terms in the sum. 
Let $f(T):={M^*(T)\over T^{1/2}}$, and 
$g_i(T):={N_i^{**}(T)\over T^{1/2}}$. We obtain the following 
pictures

\begin{figure}[H]  
\centering
\includegraphics[trim = 0mm 20mm 0mm 15mm, clip, scale=0.4]{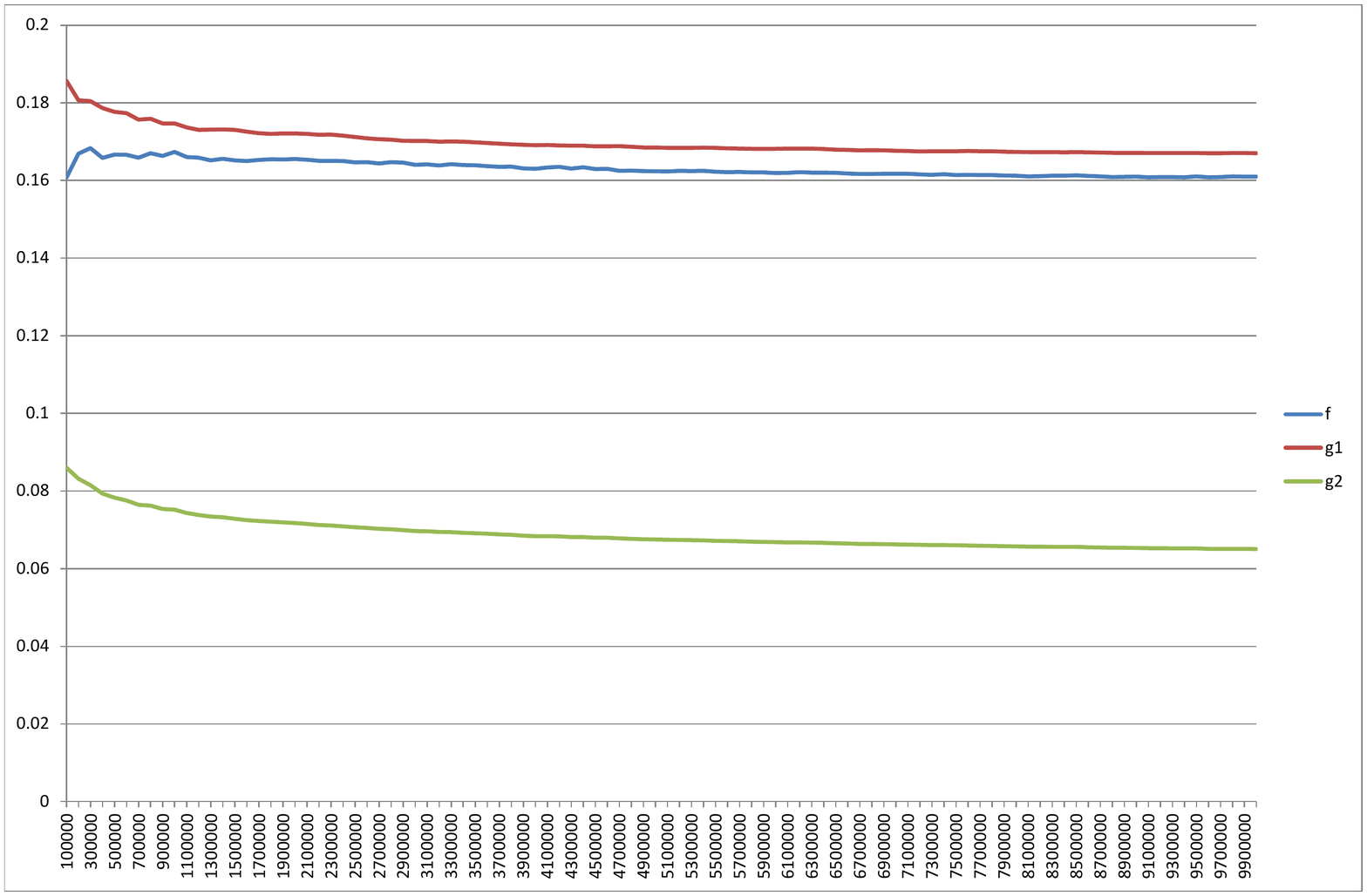} 
\caption{Graphs of the functions $f(T)$ and $g_i(T)$, $i=1,2$.} 
\end{figure}

Note similarity with the predictions by Delaunay \cite{Del0} 
for the case of quadratic twists of a given elliptic curve (and 
numerical evidence in \cite{djs} \cite{DSz}).

\subsection{The rank one case}

Let $T(X):={2\over X^*}\sum {L'(E_1(u),1)\over \Omega_{E_1(u)}}$, 
where the sum is over integers  $u\equiv 1 (\text{mod}\, 4)$, 
$|u|\leq X$, such that $u^2+64=p_1 \cdots p_k$ is a product of even 
number of different primes, 
and $X^*$ denotes the number of terms in the sum. 
Let $u(X):={T(X)\over X^{1/2}\log(X)}$. Then, using PARI/GP for 
computations of $L'(E_1(u),1)$, we obtain the following 
picture 

\begin{figure}[H]  
\centering
\includegraphics[trim = 0mm 20mm 0mm 15mm, clip, scale=0.4]{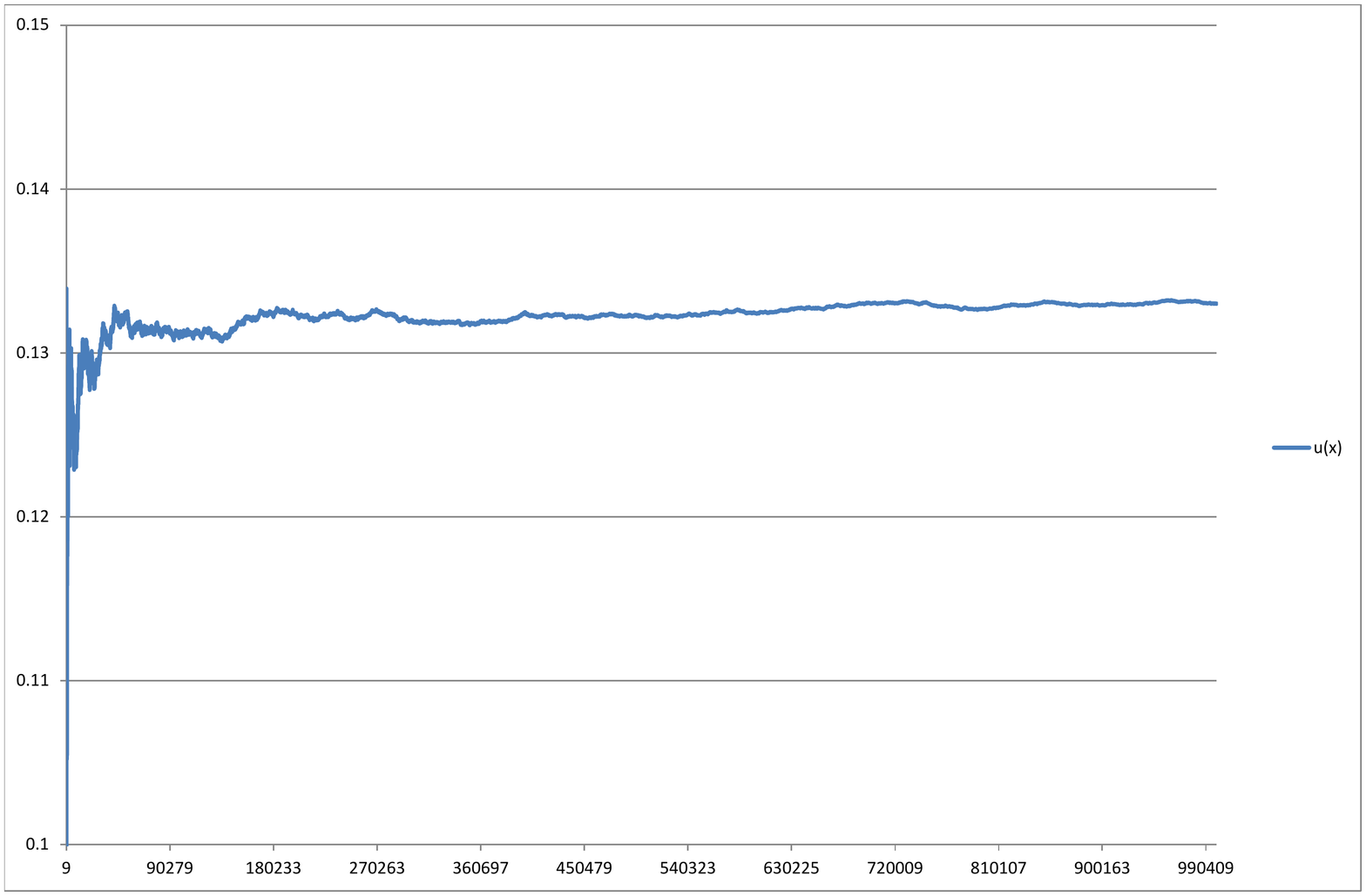} 
\caption{Graph of the function $u(X)$.} 
\end{figure} 

Hence, assuming the exact Birch and Swinnerton-Dyer conjecture 
for the rank one families $E_i(u)$, $i=1,2$, where 
$u^2+64=p_1 \cdots p_k$ is a product of even number of different 
primes, we expect the asymptotic formulae 
$$
{1\over X^*}\sum |\sza(E_i(u))| R(E_i(u)) 
\sim c_iX^{1/2}\log X, \quad \text{as} \quad X\to\infty, 
$$
where we sum over $|u|\leq X$, $u\equiv 1 (\text{mod} \,4)$, 
such that $u^2+64=p_1 \cdots p_k$ is a product of even 
number of different primes (compare \cite{djs}, section 7.2).

\bigskip 
\noindent 
{\bf Remark.} 
Delaunay and Roblot \cite{DelR} investigated regulators of 
elliptic curves with rank one in some families of quadratic 
twists of a fixed elliptic curve, and formulated some conjectures 
on the average size of these regulators. Delaunay asked us 
to do similar calculations for our family $E_i(u)$.  
We hope to consider such investigations in some future.

\section{Distributions of $L(E(u),1)$ and $|\sza(E(u))|$}

\subsection{Distribution of $L(E(u),1)$} 

It is a classical result (due to Selberg) that the values of 
$\log |\zeta({1\over 2}+it)|$ follow a normal distribution.

Let $E$ be any elliptic curve defined over $\Bbb Q$. 
Let $\cal E$ denote the set of all fundamental discriminants $d$ with $(d,2N_E)=1$ and 
$\epsilon_E(d)=\epsilon_E\chi_d(-N_E)=1$, where $\epsilon_E$ is the root number of $E$ 
and $\chi_d=(d/\cdot)$.  Keating and Snaith \cite{KS2} have conjectured that, for 
$d\in\cal E$, the quantity $\log L(E_d,1)$ has a normal distribution with mean 
$-{1\over 2}\log \log|d|$ and variance $\log \log|d|$; see \cite{CKRS} \cite{djs} 
\cite{DSz} for numerical data towards this conjecture.  

Below we consider the family of Neumann-Setzer type elliptic 
curves. Our data suggest that the 
values $\log L(E(u),1)$ also follow an approximate normal  distribution. 
Let $B=10^7$, 
$W = \{|u| \leq B : u\equiv 1 (\text{mod}\, 4) \, 
\text{and  satisfies (**)} \}$
and $I_x = [x,x+0.1)$ for $x\in\{ -10, -9.9, -9.8, \ldots, 10 \}$.
We create a histogram with bins $I_x$ from the data 
$\left\{ \left(\log L(E(u),1) + \frac 1 2 \log\log |u|\right)/\sqrt{\log\log |u|} : |u|\in W\right\}$.
Below we picture this histogram. 

\begin{figure}[H]  
\centering
\includegraphics[trim = 0mm 20mm 0mm 15mm, clip, scale=0.4]{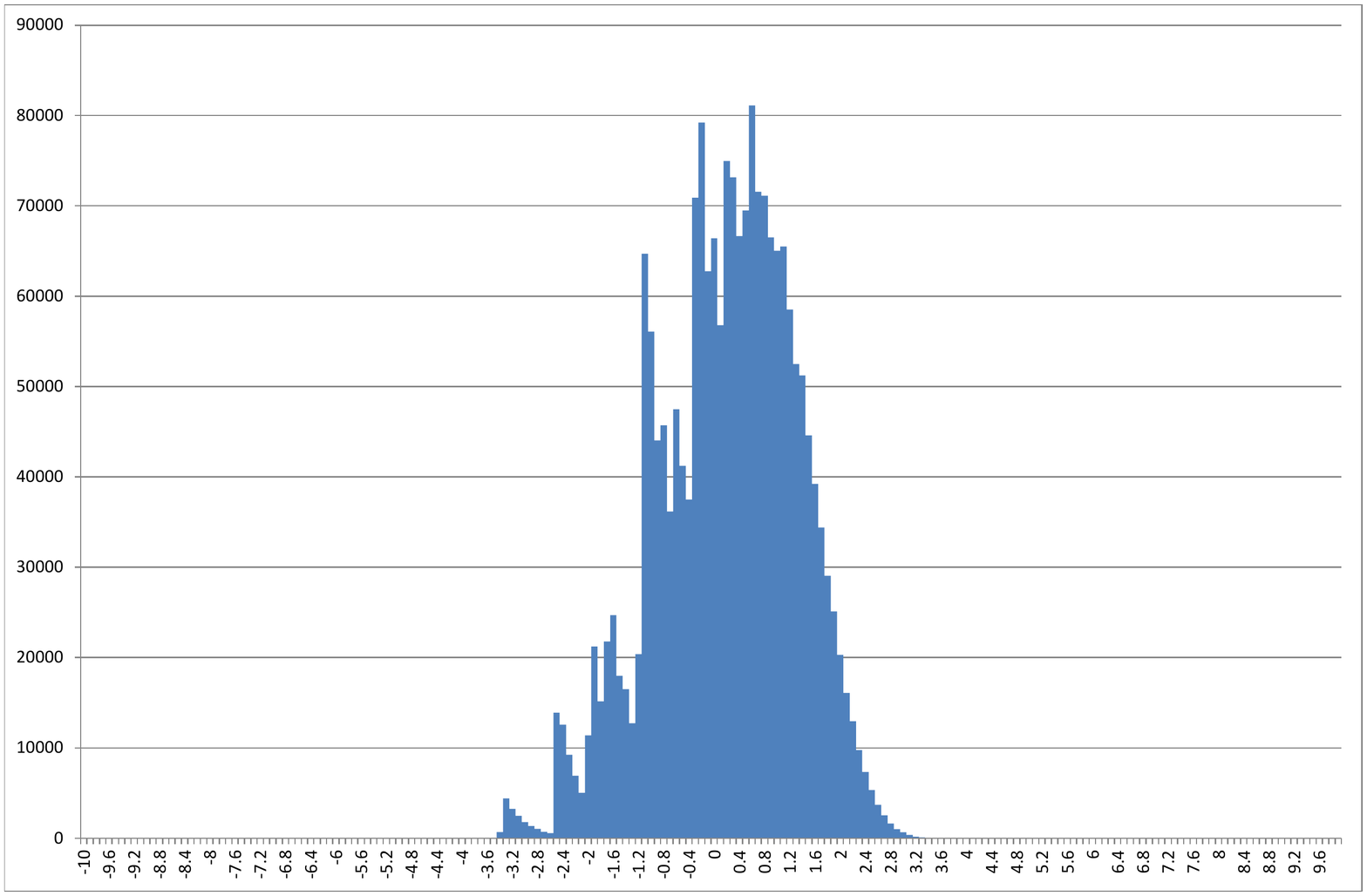} 
\caption{Histogram of values $\left( \log L(E(u),1) + {1\over 2} 
\log\log |u| \right)/\sqrt{\log\log |u|}$ for $|u| \leq B$ : $u\equiv 1 (\text{mod}\, 4)$ satisfying (**), and such that $L(E,1)\not=0$.} 
\end{figure}

\subsection{Distribution of $|\sza(E(u))|$} 

It is an interesting question to find 
results (or at least a conjecture) on distribution of the order of the Tate-Shafarevich 
group for rank zero Neumann-Setzer type elliptic curves 
$E_1(u)$ and $E_2(u)$. It turns out that the values of $\log(|\sza(E_i(u))|/\sqrt{|u|})$ 
are the natural ones to consider (compare Conjecture 1 in  \cite{RS}, 
and numerical experiments in \cite{djs} \cite{DSz}). 
Below we create histograms from the data 
$\left\{ \left( \log(|\sza(E_i(u))|/\sqrt{|u|}) - \mu_i  
\log\log |u| \right)/\sqrt{\sigma_i^2\log\log |u|} : \, |u|\in W\right\}$, 
where $\mu_1=-{1\over 2}$, $\mu_2=-{1\over 2}-\log 2$, $\sigma_1^2=1$, and 
$\sigma_2^2=1+(\log 2)^2$ (here we use Lemma 1(iii) above, and Lemma 4 in 
\cite{RS}). Our data suggest that the 
values $\log(|\sza(E_i(u))|/\sqrt{|u|})$ also follow an approximate normal 
distribution. Below we picture these histograms.

\begin{figure}[H]  
\centering
\includegraphics[trim = 0mm 20mm 0mm 15mm, clip, scale=0.4]{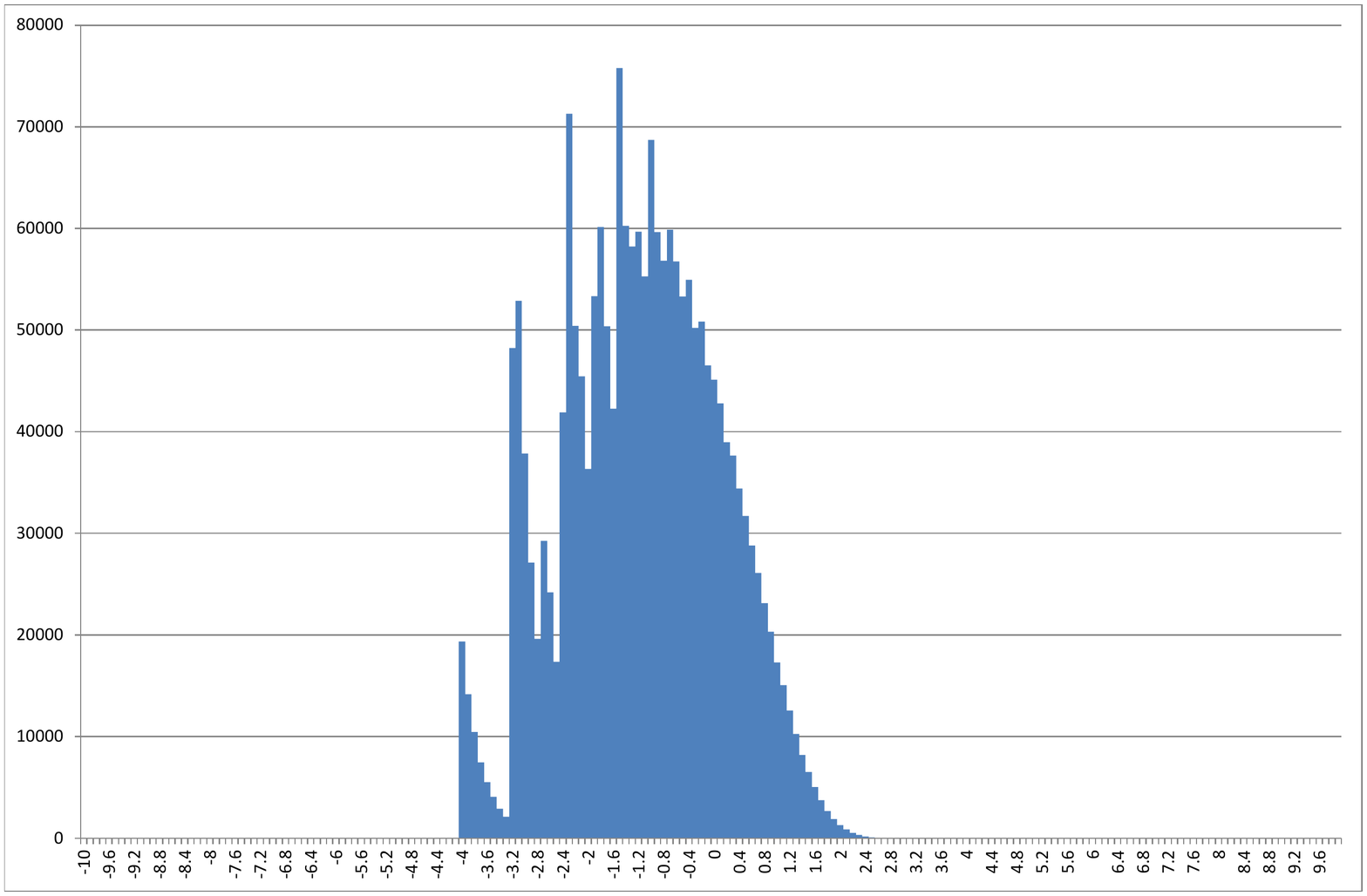} 
\caption{Histogram of values $\left( \log(|\sza(E_1(u))|/\sqrt{|u|}) + {1\over 2} 
\log\log |u| \right)/\sqrt{\log\log |u|}$ for $|u| \leq B$ : $u\equiv 1 (\text{mod}\, 4)$ satisfying (**), and such that $L(E,1)\not=0$.} 
\end{figure}

\begin{figure}[H]  
\centering
\includegraphics[trim = 0mm 20mm 0mm 15mm, clip, scale=0.4]{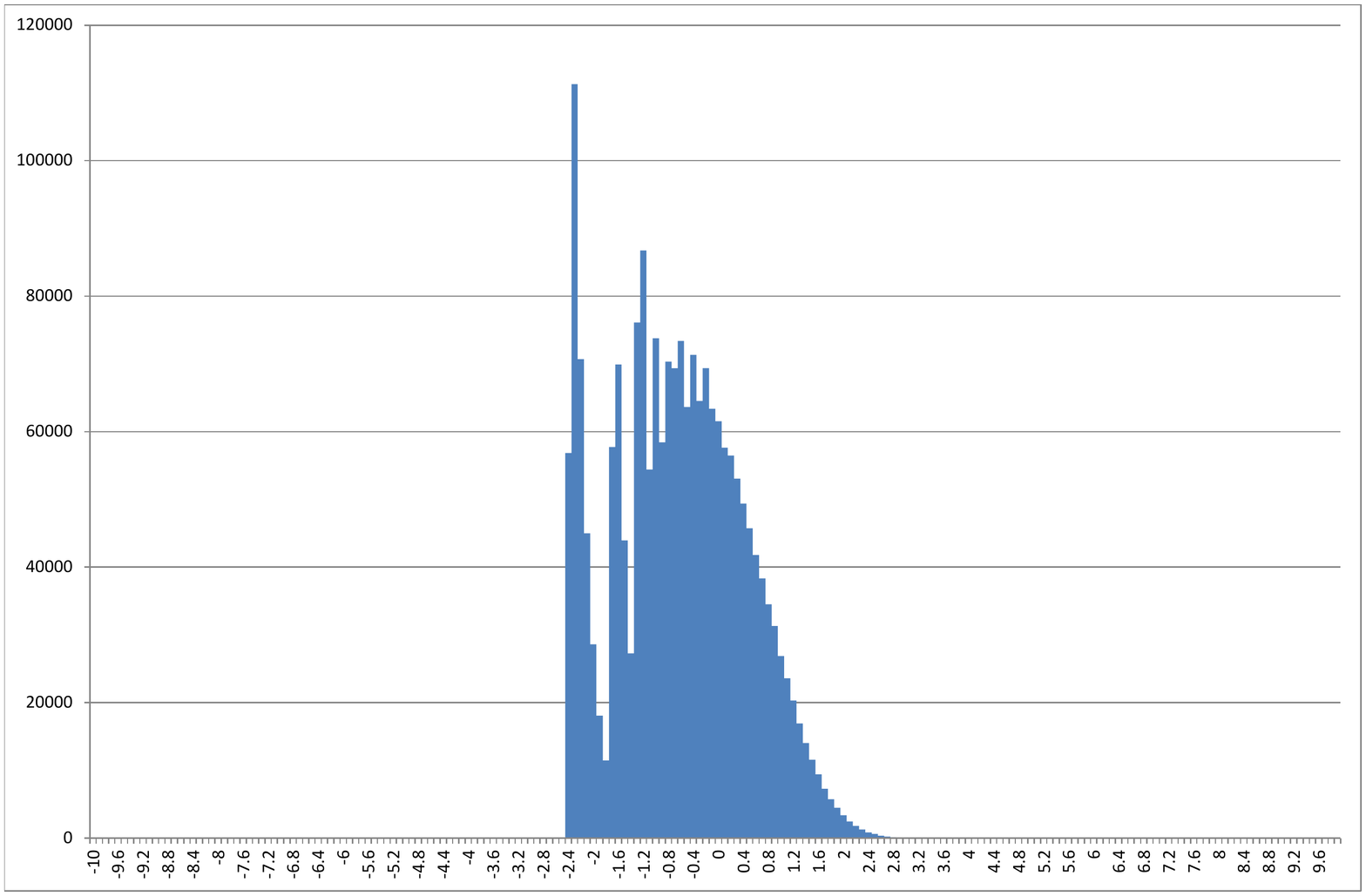} 
\caption{Histogram of values $\left( \log(|\sza(E_2(u))|/\sqrt{|u|}) + ({1\over 2}+\log 2)  
\log\log |u| \right)/\sqrt{(1+(\log 2)^2)\log\log |u|}$ for $|u| \leq B$ : $u\equiv 1 (\text{mod}\, 4)$ satisfying (**), and such that $L(E,1)\not=0$.} 
\end{figure}

Institute of Mathematics, University of Szczecin, Wielkopolska 15, 
70-451 Szczecin, Poland; E-mail addresses: dabrowsk@wmf.univ.szczecin.pl and dabrowskiandrzej7@gmail.com;  
lucjansz@gmail.com

\end{document}